\documentclass[a4paper, 11pt]{article}
\usepackage{amsmath,amsthm,amsfonts,amssymb}
\usepackage{ifthen}

 \provideboolean{showlabels}
 \setboolean{showlabels}{false}
\newcommand{\mylabel}[1]{\ifthenelse{\boolean{showlabels}}{{\tt{[{#1}]}}\label
{#1}}{\label{#1}}}

\theoremstyle{plain}

\newtheorem{thm}{Theorem}[section]
\newtheorem{lem}[thm]{Lemma}
\newtheorem{lemma}[thm]{Lemma}
\newtheorem{pro}[thm]{Proposition}
\newtheorem{cor}[thm]{Corollary}
\newtheorem*{A}{Theorem A}

\theoremstyle{definition}
\newtheorem{rem}{Remark}

\newtheorem*{axiom}{Inductive Hypothesis}
 \newtheorem{claim}{Claim}

\newcommand{\F}{{\mathcal{F}}}
\newcommand{\B}{{\mathcal{B}}}
\newcommand{\RB}{{\mathcal{R}}}
\newcommand{\SB}{{\mathcal{S}}}
\newcommand{\UB}{{\mathcal{U}}}
\newcommand{\VB}{{\mathcal{V}}}
\newcommand{\WB}{{\mathcal{W}}}
\newcommand{\XB}{{\mathcal{Y}}}
\newcommand{\Xii}{{\mathcal{U}}}
\newcommand{\E}{{\mathcal{E}}}

\newcommand{\TB}{{\mathcal{T}}}
\newcommand{\ZZ}{\mathbb{Z}}
\newcommand{\W}{{\mathcal{V}}}
\newcommand{\D}{{\mathcal{D}}}

\newcommand{\beq}[1]{\begin{equation}\mylabel{#1}}
\newcommand{\eeq}{\end{equation}}
\newcommand{\comment}[1]{}

\DeclareMathOperator{\Irr}{Irr}

\DeclareMathOperator{\Gal}{Gal}
\DeclareMathOperator{\Hom}{Hom}

\DeclareMathOperator{\End}{End}

\begin{document}
\centerline{\bf{Hyperbolic modules and cyclic subgroups}}
\bigskip

\centerline{\bf{Maria Loukaki}}
\bigskip

\centerline{Dept. of Applied  Mathematics, University of Crete,}
\centerline{ Knossos Av. GR-71409, Heraklion-Crete, Greece}
\bigskip

\centerline{E-mail: loukaki@tem.uoc.gr }

\abstract
Let $G$ be a finite group of odd order, $\F$ a finite field of odd characteristic $p$
 and $\B$ a  finite--dimensional     symplectic $\F G$-module.
We show that $\B$ is $\F G$-hyperbolic, i.e., it contains a self--perpendicular $\F G$-submodule,
iff it is $\F N$-hyperbolic for every cyclic subgroup $N$ of $G$.

\section{Introduction}
Let $\F$ be a finite field of odd characteristic $p$, $G$ a finite group  and $\B $ a finite--dimensional  $\F G$-module.
If $\B $ carries a non-singular  alternating bilinear form $<\cdot, \cdot>$ (i.e.,  a symplectic form)
   that is invariant   by $G$, then we call $\B$ a {\em symplectic } $\F G$-module.
Following the notation in \cite{1},  for any   $\F G$-submodule  $\SB$ of $\B$,  we write  $\SB ^{\perp}$  for 
the { \em perpendicular  }  subspace  of $\SB$, i.e.,  $\SB^{\perp}:= \{t \in \B |<\SB , t> = 0\}$. 
We say that $\SB$ is {\em isotropic }  if $\SB \leq \SB^{\perp}$, and $\B $ is   { \em  anisotropic } 
if it contains no non--trivial isotropic 
$\F G$-submodules. Furthermore,  we say that $\B$ is  { \em hyperbolic } if it contains some
 self--perpendicular $\F G $-submodule $\SB$, i.e., $\SB$ is an
 $\F G$-submodule satisfying  $\SB = \SB^{\perp}$.  

Symplectic modules play an essential role in studying monomial characters. (An irreducible character
$\chi$ of a finite group $G$  is monomial if it is induced from a linear character of a subgroup of $G$.)
One of the most representative links between symplectic modules and monomial characters  
can be found in \cite{1}. (For other examples one could look at   \cite{9,4, 5, 6, 11, 10,  7}, and 
\cite{8}.) 
There E. C. Dade proved the following theorem   (Theorem 3.2 in \cite{1}): 
\begin{thm}[Dade]
Suppose  that $\F$ is a finite field of odd characteristic $p$, that $G$ is a finite $p$-solvable group, that $H$
is a subgroup of $p$-power  index  in $G$, and that $\B$ is a symplectic $\F G$-module whose restriction  $\B _H$ 
to a symplectic $\F H$-module is hyperbolic. Then $\B$ is hyperbolic.
\end{thm}
Using the above theorem, E. C. Dade   was able to prove  (Theorem 0 in \cite{1}) that,    given a $p$-solvable odd group $G$, 
 an irreducible monomial  character $\chi$ of $G$, and a subnormal subgroup  $N$ of $G$,  
every irreducible constituent of the restricted character $\chi_N$
is monomial, provided that $\chi(1)$ is a power of $p$.

In this paper we prove 
\begin{A}
Suppose that $\mathcal{F}$ is a finite field of odd characteristic $p$,
 that $G$ is a finite group of odd order, and that $\mathcal{B}$ is a symplectic 
${\mathcal{F}} G$--module whose restriction $ {\mathcal{B}}_{N}$   to a symplectic 
${\mathcal{F}}N$--module is hyperbolic for every cyclic subgroup $N$ of $G$.
Then $\mathcal{B}$ is hyperbolic.
\end{A}

All groups considered here are of finite order, and all modules have finite dimension 
over  $\F$.

{\bf Acknowledgments}
I am indebted to Professor E. C. Dade for  many   helpful ideas and 
suggestions.
 Also,  I would like to thank  Professor M. Isaacs for  useful  conversations 
 that helped me improve this paper.

\section{Symplectic modules }
We first give some elementary results about symplectic modules.

Assume that $\B$ is a symplectic $\F G$-module, while $\SB$ is an isotropic 
$\F G$-submodule  of $\B$.  Then the factor $\F G$-module 
$\bar{\SB}= \SB^{\perp}/ \SB$ is again a symplectic $\F G$-module  with the symplectic form 
defined as (see 1.4 in \cite{1}), 
\begin{equation}\mylabel{e10} 
< s_1 + \SB , s_2 +\SB> = <s_1, s_2>, \text{ for all } s_1, s_2 \in \SB^{\perp}. 
\end{equation}
Furthermore,  if $\SB$ is an isotropic $\F G$-submodule of $\B$, then  its  $\F$-dimension 
$\dim_{\F} \SB $ is at  most  $(1/2)  \dim_{\F} \B$, (see 19.3 in \cite{as}).

We  say  that an isotropic $\F G$-submodule  $\SB$  of $\B$ is 
{ \em maximal isotropic }  if $\SB$  is  not properly  contained in any larger 
isotropic $\F G$-submodule of $\B$.  
Clearly any self--perpendicular 
$\F  G$-submodule  $\SB$ of  $\B$ is maximal isotropic. 
The converse is also correct under the extra assumption that $\B$ is  $G$-hyperbolic (see Lemma 3.1 in \cite{1}).
Another way  to get a self--perpendicular module from a maximal isotropic one is to control its dimension, as 
the following lemma shows.
\begin{lemma}\mylabel{ll1}
Assume that $\B$ is a symplectic $\F G$-module, and that $\SB$ is a maximal isotropic 
$\F G$-submodule of $\B$.  If $\dim_{\F} \SB = (1/2) \dim_{\F}  \B$ then 
$\SB$ is self--perpendicular and $\B$ is $G$-hyperbolic.
\end{lemma}

\begin{proof}
Let $\hat{\SB}$ denote the dual of $\SB$. Then 
 $\B / \SB^{\perp} \cong \hat{\SB}$.
But $\dim_{\F} \hat{\SB}  = \dim_{\F} \SB  = (1/2) \dim_{\F} \B$.
Hence $\dim_{\F} \SB^{\perp} = (1/2) \dim_{\F} \B$.
Since $\SB \leq \SB^{\perp}$  we conclude that $\SB= \SB^{\perp}$.   
Thus the lemma holds.
\end{proof}

 The following is Proposition 2.1 in \cite{1}.
\begin{pro}\mylabel{p1}
Let $G$ be a finite group  and $\B $ be  an
anisotropic symplectic $\F G$--module. Then $\B$ is an orthogonal
direct sum:
\begin{equation}\mylabel{e01}
\B = \UB _1 \dot{\perp}\, \UB _2 \dot{\perp} \ldots \dot{\perp} \, \UB _k,
\end{equation}
where $k\geq 0$ and each $\UB _i$ is a simple $\F G$--submodule of $\B $ that is also symplectic.
\end{pro}

\begin{rem}\mylabel{rem}
If $G$ has odd order then according to Proposition (1.10) and
Corollary 2.10 in \cite{1} all the $\UB _i$ that appear in (\ref{e01})
are distinct.
\end{rem}

\begin{lem} \mylabel{intr1}
Let $\UB $ be an $\F G$-module that affords a symplectic $G$-invariant  form $< \cdot , \cdot >$. 
Then $\UB$ is self--dual.
\end{lem}

\begin{proof}
We write $\widehat {\UB }$ for the dual  $\F G$-module of $\UB$. 
 For every $x \in \UB$ the map $\alpha_x :\UB \to \F$ 
defined as: 
 $$ 
\alpha _x(u) = <u\, ,\, x> \, \text{for all }  u \in \UB 
$$
 is an element of $ \Hom _{\F}(\UB ,\F) \cong \widehat{\UB }$.  Since   $<\cdot \, ,\, \cdot>$  is $G$-invariant 
the map $\alpha: x \to \alpha_x$  is an $\F G$-homomorphism   from $\UB$  to $\widehat{\UB}$.
Furthermore the kernel  of
$\alpha$ is trivial,  as $\UB $ is symplectic. Hence  $\UB \cong \widehat{\UB }$. 
\end{proof}

\begin{cor}\mylabel{co1}
Let $\B$ be an anisotropic symplectic $\F G$--module. Then each  of
the simple $\F G$--modules $\UB _i$  that appears in \eqref{e01} is self--dual.
\end{cor}
\begin{proof}
It follows easily from Proposition \ref{p1} and  Lemma \ref{intr1}. 
 \end{proof}

\begin{pro}\mylabel{intr2}
Assume that $\UB$ is a simple symplectic $\F G$--module. Let $N$ be a
normal subgroup of $G$ such that $|G:N|$ is odd. Then any simple $\F N$--submodule of 
$\UB _N$ is self--dual.  Hence any $\F N$-submodule of $\UB _N$ 
is self--dual.
\end{pro}
\begin{proof}
As $N$ is a normal subgroup of $G$, Clifford's theorem  implies that
\begin{equation}\mylabel{du}
 \UB _N  \cong e( \VB _1 \oplus
\ldots \oplus \VB _n)
\end{equation}
where $\VB =\VB _1$ is a
simple  $\F N$--submodule of $\UB$ and $\VB _1, \ldots , \VB _n$ are
the distinct $G$-conjugates of $\VB$. 
 So $n \big{|} |G:N| $ and therefore
$n$ is odd.

According to Lemma
\ref{intr1} the module $\UB$ is self--dual. Hence the dual, $\widehat{\VB _i}$ of
any $\VB _i$ should appear in (\ref{du}). Therefore we can form pairs
among the  $\VB _i$, consisting of a simple $\F
N$--module $\VB _k$ and its dual for $k \in \{ 1, \ldots , n \}$,
where we take as the second part of the pair the module itself if it
is self--dual.
Since $G$ acts transitively on  the $\VB _i$ for $i= 1, \ldots , n$, 
either all the $\VB _i$ are self--dual or none of them is. In the
latter case we get that any of the above pairs consists of two
distinct modules. This implies that  $2 | n$.
As $n$ is odd, this case can never occur. Hence any one of the $\VB_i$ is 
self--dual and the proposition is proved.
\end{proof}

\begin{pro}\mylabel{intr3}
Assume that the symplectic $\F G$--module $\B$ is hyperbolic. Assume further that 
$\B$ is a semi--simple  $\F G$-module. Then
every self--dual  simple $\F G$--submodule of $\B$ appears with even multiplicity in any
decomposition of $\B$ as a direct sum of simple $\F G$--submodules.
\end{pro}
\begin{proof}
Because $\B$  is hyperbolic  it contains a self--perpendicular 
$\F G$-submodule $\SB$.    For every $\F G$-submodule $\VB$ of $\B$ we have 
$\B /  \VB^{\perp}  \cong \widehat{\VB}$. 
So 
\begin{equation}\mylabel{e11}
\B / \SB  \cong \widehat{\SB}
\end{equation}
Now the proposition  follows from \eqref{e11}  and the fact that $\B$ is semi--simple.
\end{proof}

\begin{cor}\mylabel{intr4}
Let  $\B$ be  an anisotropic  symplectic 
$\F G$--module. Let $N$ be a
normal subgroup of $G$ such that $|G:N|$ is odd.
Assume further that $\B_N$ is a hyperbolic $\F N$-module. 
 Then any simple $\F N$--submodule of $\B _N$ appears with even multiplicity in any 
decomposition of  $\B _N$ as a direct sum of simple $\F N$--submodules.
\end{cor}
\begin{proof}
This  is a straightforward application of Propositions \ref{p1},  \ref{intr2} and \ref{intr3}.
\end{proof}

We close this section with a well known fact that we prove here for completeness.
\begin{lem}\mylabel{intr5}
Assume that $\UB$ is a self--dual absolutely irreducible $\F G$-module, where  $G$ has odd order and  
$\F$ is a finite field whose characteristic does not divide $|G|$.
Then $\UB$ is trivial.
\end{lem}

\begin{proof}
Let $\chi$  denote the $\F$-absolutely irreducible character that $\UB$ affords, while  $\phi $ denotes  a  Brauer character 
that  $\UB$ affords.  Then $\phi$ is defined for every element of $G$, since the characteristic of $\F$ is coprime to $|G|$.
Because $\UB$ is self--dual, the character $\phi $ is real valued.  
Let $\nu_2(\phi) = |G|^{-1} \sum_{g \in G}\phi(g^2)$ be the Frobenius--Schur indicator (see Chapter 4 in \cite{3}) of $\phi$.
Then  Theorem 4.5 in \cite{3}  implies that $\nu_2(\phi )\ne  0$, since $\phi$ is real valued.
But 
$$
\nu_2(\phi) =  |G|^{-1} \sum_{g \in G}\phi(g^2) =  |G|^{-1} \sum_{g \in G}\phi(g), 
$$
because $G$ has odd order.
Hence $\nu_2(\phi)$ is the inner product $\nu_2(\phi)  = [ \phi , 1_G]$, where $1_G$ is the trivial character of $G$.
We conclude that $[\phi, 1_G] \ne 0$. Hence $\phi = 1_G$. Therefore $\chi= 1_G$,  and the 
 lemma follows. 
\end{proof}

\section{ Proof of Theorem A }
We can now prove our main result.  
The proof will follow from a series of lemmas, based on the
hypothesis that  ${\mathcal{F,B}}, G$ form a minimal
counter--example. All the groups considered in this section have 
odd order. We also fix  the odd prime $p$ that is the characteristic of $\F$,  and we assume that

\begin{axiom} 
 {\sl ${\mathcal{F,B}}, G$ have been
chosen among all triplets satisfying the hypothesis but not the
conclusion of  Theorem A so as to minimize first the order $|G|$ 
 of $G$ and then the  $\mathcal{F}$--dimension $\dim_{\F}\mathcal{B}$ of
$\mathcal{B}$.}
\end{axiom}

\begin{rem} \mylabel{r1}
For any proper subgroup $H$ of $G$ the minimality of $|G|$ 
easily implies that the restriction   ${\mathcal{B}}_{H} $ is  a hyperbolic
$\F H$-module.
\end{rem}

\begin{lem} \mylabel{l1}
 $\mathcal{B}$ is non--zero and anisotropic.
\end{lem}

\begin{proof}
If  $\mathcal{B}$ were zero it would be hyperbolic contradicting
 the Inductive Hypothesis. So $\mathcal{B}$ is non--zero.
If   $\mathcal{B}$ is not anisotropic then it contains a non--zero
 isotropic $\mathcal{F}$$ G$--module  $\mathcal{U}$. Let $N$ be an arbitrary  cyclic subgroup of $G$.
Then the isotropic $\F N$-submodule $\mathcal{U}_N$  of $\B _N$ 
is contained in  some maximal isotropic
 ${\mathcal{F}}N$--submodule  $\mathcal{V}$ of ${\mathcal{B}}_{N}$. 
  Since  ${\mathcal{B}}_{N}$ is hyperbolic this maximal isotropic submodule 
 is self--perpendicular, i.e., 
   $\mathcal{V}=\mathcal{V}^{\perp} $. Hence 
 $$
\mathcal{U} \subseteq \mathcal{V}=\mathcal{V}^{\perp}  \subseteq \mathcal{U}^{\perp}. 
$$
Therefore the factor module
  $\bar{\mathcal{V} } ={\mathcal{V}}/{\mathcal{U} } $
 is a self--perpendicular ${\mathcal{F}}N$--submodule of
 the symplectic ${\mathcal{F}}G$--module 
$\bar{\mathcal{U}} = \mathcal{U}^{\perp} /{\mathcal{U}}$.
Hence  ${\mathcal{F}}, G, \bar{\mathcal{U}} $ 
satisfy the hypothesis of the Main Theorem. 
As $\dim( \bar{\mathcal{U} }) <  \dim(\mathcal{B}) $, the minimality of
$ \dim(\mathcal{B}) $ implies that $\bar{\mathcal{U}}$ is a hyperbolic 
 ${\mathcal{F}}G$--module. So there is a self--perpendicular 
${\mathcal{F}}G$--submodule $\bar{\mathcal{J}}$ in $\bar{\mathcal{U}}$. 
From the definition of the symplectic form on $\bar{\mathcal{U}}$ (see \eqref{e10})
it follows that the inverse image $\mathcal{J} $ of $\bar{\mathcal{J}}$
 in $\mathcal{U}^{\perp} $ is a self--perpendicular
 ${\mathcal{F}}G$--submodule of  $\mathcal{B}$ containing
 $\mathcal{U}$. Therefore  $\mathcal{B}$ is
hyperbolic, contradicting the Inductive Hypothesis.
So the lemma holds.
\end{proof}

\begin{lem} \mylabel{l2}
$p$ doesn't divide the order $|G|$ of $G$.
\end{lem}

\begin{proof}
Suppose  that $p $ divides $|G| $.  Because  $G$  is solvable,  
it contains  a Hall $p'$--subgroup $H$.  If $G$ is a $p$-group we take $H = 1$.
 Since  $p$ divides $|G|$, the  subgroup $H$ is  strictly smaller than $G$.
  Then according to Remark \ref{r1},   the $\F H$-module $\mathcal{B}$$_{H}$ is
hyperbolic. It follows (see Theorem 3.2 of \cite{1} ) that $\mathcal{B}$ is
a hyperbolic $\F G$-module, contradicting the Inductive Hypothesis. 
Hence $(p, |G|) =1$.
\end{proof}

\begin{lem}\mylabel{l3}
$\mathcal{B}$ is an orthogonal direct sum
\begin{equation}\mylabel{e-1}
{\mathcal{B}} = {\mathcal{U}}_{1} \dot{\perp}  \ldots \dot{\perp}\,{\mathcal{U}}_{k}
\end{equation}
where $k\geq 1$ and  $\bigl\{ {\mathcal{U}}_i\bigr\}_{i=1, \ldots ,k}$ are
distinct, simple ${\mathcal{F}}G$--submodules of
$\mathcal{B}$, that are also  symplectic.
Furthermore each  ${ \mathcal{U}}_i$ is quasi--primitive (i.e.,  its
restriction to every normal subgroup of $G$ is homogeneous).
\end{lem}

\begin{proof}
The first statement follows from Lemma~\ref{l1},  Proposition \ref{p1}  and Remark \ref{rem}.
 For the rest of the proof we fix 
${ \mathcal{U}}={\mathcal{U}}_i $ 
for some $ i=1, \ldots , k$.
We also fix a normal
subgroup $K$ of $G$. If the restriction of $ \mathcal{U}$ to  $K$ is
not homogeneous then Clifford's Theorem implies that 
$$
{\mathcal{U}}_K \cong e( {\mathcal{V}}^{\sigma_1} \oplus
\ldots \oplus {\mathcal{V}}^{\sigma_r} )
$$ 
 where $e$ is some positive integer,  ${\mathcal{V}} ={\mathcal{V}}^{\sigma_1}$ is a
 simple $ {\mathcal{F}}K$--submodule of
$\mathcal{U}$ and ${\mathcal{V}}^{\sigma_1}, \ldots
,{\mathcal{V}}^{\sigma_r} $ are the distinct conjugates of
$\mathcal{V}$ in $G$, with ${\sigma_1 , \ldots , \sigma_r}$ coset
representatives of the stabilizer, $G_{\mathcal{V}}$,  of $\mathcal{V}$
in $G$.

Let $\WB  =\UB (\VB) $ be the $\VB$-primary component of $\UB_K$.
 Then Clifford's Theorem implies that $\WB$ is
the unique simple $\F G _{\VB}$--submodule of $\UB$
that lies above $\VB$ and induces $\UB$,  i.e., that satisfies
 $\WB ^G \cong \UB$ and
$\WB _K \cong e\VB $.
Furthermore the dual $\widehat{\WB} $ of $\WB$ induces in $G$ the dual
$\widehat{\UB}$ of $\UB$ since  $\widehat{\WB} ^ G \cong \widehat{\WB ^G}$. 
Hence $\widehat{\WB} ^G   \cong    \UB$,  because 
$\UB $ is self--dual (see Lemma \ref{intr1}).
On the other hand, the restriction of $\widehat{\WB}$ to $K$ is isomorphic
to $e\VB$ since  $\VB$ is self--dual by Proposition \ref{intr2}. 
Hence  the unicity of $\WB$ implies that $\WB$  is self--dual.

According to  Proposition \ref{intr3} the self--dual $\F G_{\VB}$-module 
$\WB $ appears with even multiplicity as a direct summand of
$\B _{G_{\VB}}$, because  $\B _{G_{\VB}}$ is hyperbolic ($G_{\VB} < G$).
 This, along with the fact that 
$\WB$ appears with multiplicity one in $\UB _{G_{\VB}}$, implies that
 there is some $j \in  \{ 1, \ldots , k\} $
 with $j\ne i$ such that the $\VB$-primary component $\UB _j(\VB) $ of $\UB_j$ 
 is  isomorphic to $\WB$.  So 
$$
\WB =  \UB(\VB) \cong \UB _j(\VB).
$$ 
We conclude that 
$$
\UB_i  = \UB \cong \WB ^G \cong  \UB _j(\VB)^G \cong \UB _j, 
 $$
as $\F G $-modules. 
 This contradicts the fact that  $\{\UB_i \}_{i=1}^k$  are all distinct, by the first statement of the lemma.
Hence the lemma is proved.
\end{proof}

From now on and until the end of the paper,  we write   $\E$ for a finite algebraic  field extension of $\F$, 
that is a splitting field of $G$ and all its subgroups.

\begin{lem}\mylabel{l3.3}
Assume that $ \UB_i$, for $i=1, \dots, k$, is a direct summand of $\B$  appearing in \eqref{e-1}.
Let $N \unlhd G$. Then $\UB_i |_N  \cong  e_i \VB_i $, where $\VB_i$ is an irreducible 
$\F N$-submodule of $\UB_i$  and $e_i$ is an  integer. If $\VB_i$ is non-trivial 
then $e_i$ is odd.  
\end{lem}

\begin{proof}
We fix  $\UB = \UB_i$, for some $i=1, \dots, k$. We also fix  a normal subgroup $N$ of $G$.
According to Lemma \ref{l3}, the $\F G$-module $\UB$ is quasi--primitive.
Hence there exists an irreducible $\F N$-submodule $\VB$  of $\UB$, and an integer $e$  such that 
$\UB _N \cong e \VB$. Thus, it remains to show that $e$ is odd in the case that $\VB$ is  non-trivial.
So we assume that $\VB$,  and thus $\UB$,   is   non-trivial. 

We observe that if $\UB$ and $\VB$  were  absolutely irreducible modules  then 
it would be  immediate  that $e$ is odd (even  if $\VB$ was  trivial),  because   for absolutely irreducible modules the integer 
$e$ divides the order of $G$ (see Corollary 11.29 in \cite{3}). So we assume that $\F$ is not a splitting field of $G$, and we work 
with the algebraic  field extension $\E$ of $\F$.
We define $\UB^{\E}$   to be the extended $\E G$-module 
$$
\UB^{\E} = \UB \otimes_{\F} \E. 
$$
According to Theorem 9.21 in \cite{3}, there exist absolutely irreducible $\E G$-modules $\UB ^i $, for $i=1, \dots, n $,   such that 
\begin{equation*}\mylabel{ee1}
\UB ^{\E} \cong \bigoplus_{i =1}^n \UB ^i .
\end{equation*}
Furthermore  the $\UB ^i $, for all  $i=1, \dots, n$,  constitute a Galois conjugacy class over $\F$, and thus  they are all distinct.
In particular,  if $\E _{\UB} $ is the 
subfield of $\E $  that is generated by  all the values of the irreducible $\E$-character  afforded by $\UB^i$  
(the same  field for all $i=1, \dots, n$),
then $n = [\E_{\UB} : \F ] = \dim_{\F} (\E_{\UB})$.
(Note that  $\E_{\UB} $ is  the unique subfield of $\E$ isomorphic 
to the center of the endomorphism algebra $\End _{\F G }(\UB)$.)  
Clearly $n \cdot \dim_{\E} \UB ^1 = \dim_{\F} \UB$.  Hence $n$ is even,  because   $\dim_{\F} \UB $ is even (as $\UB$ is symplectic)
and $\dim_{\E} \UB^1$ is odd  (as $G$ is odd and $\UB^1$ is an absolutely irreducible 
$\E G$-module).
In addition, each  $\E G$-module $\UB^i$,  for $i=1, \dots, n$, when consider as an $\F G$-module, 
is isomorphic  to a  direct sum of $[ \E : \E _{\UB}]$ copies of  $\UB$ (see Theorem 1.16  in Chapter 1 of  \cite{hu}). 
Hence if we denote by  $\UB^i _{\F }$ the $\E G$-module  $\UB^i$ regarded as an  $\F G$-module,  we get 
\begin{equation}\mylabel{a}
\UB^i _{\F} \cong  [\E : \E_{\UB}]  \UB ^i, 
\end{equation}
for all $i=1, \dots, n$.

We also write $\VB ^{\E}$ for the extended $\E N$-module $\VB^{ \E} = \VB \otimes_{\F} \E$.
 Then according to  Theorem 9.21 in \cite{3}  there exist absolutely irreducible $\E N$-modules $\VB ^j$ for $j=1, \dots, d$,  such that 
\begin{equation}\mylabel{ee2}
\VB^{ \E} \cong  \bigoplus_{j=1}^d \VB ^j .
\end{equation}
In addition,  the absolutely irreducible  modules $\VB ^j$,  for  all $j=1, \dots, d$, form a Galois conjugacy class, 
and thus they are all distinct.
Furthermore,  $d= [\E_{\VB} : \F] = \dim _{\F} \E_{\VB}$, where $\E _{\VB}$ is the  subfield of $\E$ generated by all the values 
of the irreducible $\E$-character afforded by $\VB^j$ (the same field  for all  $j=1, \dots, d$). 
The field $\E _{\VB}$ is the unique subfield of $\E$ isomorphic to 
the center of the endomorphism algebra $\End _{\F N}(\VB)$.
Note that, according to Proposition \ref{intr2}, the $\F N$-submodule $\VB$ of $\UB$ is self--dual.
Hence  $\VB  ^{\E} $ is also a self--dual $\E N$-module.
Because $\VB $ is non-trivial, $\VB ^j$  is also non-trivial,  for all $j=1, \dots, d$.
Therefore  the absolutely irreducible 
$\E N$-module $\VB^j$  can't be self--dual,   because $N$ has odd order and $\VB ^j$ is non--trivial (see Lemma \ref{intr5}), for all such $j$.
The fact that none of the  $\VB^j$ is self--dual, for all $j=1, \dots, d$, while they all appear in \eqref{ee2} in dual pairs, implies 
that  $d$ is  even.
Even more, if $\VB^ j_{\F}$ denotes the module $\VB^j $  regarded as an $\F N$-module, then Theorem 1.16 of Chapter 1 in \cite{hu}
implies that 
\begin{equation}\mylabel{b}
\VB^j _{\F} \cong [\E :\E_{\VB}] \VB^j, 
\end{equation}
for all $j=1, \dots, d$.

Without loss we may assume that $\VB^1, \dots , \VB^c$ are exactly 
those among the $\VB^j$, for $j=1, \dots, d$, that  lie under $\UB^1$,  for some $c=1, \dots, d$.
Thus Clifford's theorem implies that 
\begin{equation}\mylabel{ee4}
\UB^1_N   \cong e' (\VB^1 \oplus \dots \oplus \VB^c), 
\end{equation}
where $\VB^1, \dots, \VB^c$  are the distinct $G$-conjugates of $\VB^1$,   and $e', c$ are integers that divide $|G|$ and thus 
are odd. (Note that here we are dealing with  absolutely irreducible modules so $e'$ does divide $|G|$.)
If we regard the modules of \eqref{ee4} as modules over  the field $\F$ then we clearly have 
$\UB^1_{\F} |_N \cong e' (\VB^1_{\F} \oplus \dots \oplus \VB^c_{\F})$. This, along with 
\eqref{a} and \eqref{b}, implies
\begin{equation*}
[\E: \E_{\UB} ] \UB_N \cong e' c [\E:\E_{\VB}] \VB.
\end{equation*}
Since $\UB_N \cong e \VB$, we have 
\begin{equation}\mylabel{c}
[\E: \E_{\UB} ] e =  e' c [\E:\E_{\VB}]. 
\end{equation}
If $\mathcal{D}$  is the subfield of $\E$ generated by $\E_{\VB}$ and $\E_{\UB}$, then dividing both sides of 
\eqref{c} by $[\E  :\D]$ we obtain
\begin{equation}\mylabel{d}
e [\D: \E_{\UB}] = e' c [\D: \E_{\VB}].
\end{equation}

Assume that $e$ is even.  Then \eqref{d} implies that  $[\D:\E_{\VB}]$ is even, as $e'$ and $c$  are  known to be odd.
 Let $\Gamma$ be the Galois group $\Gamma = \Gal(\D / \F )$ 
of $\D$ over $\F$.  Because $\Gamma $  is cyclic, it contains a unique involution $\iota$.
Let $\E_{\VB}^*$  and $\E_{\UB}^*$ be  the subgroups of $\Gamma$ consisting of those elements of $\Gamma $ that fix 
pointwise $\E_{\VB}$ and $\E_{\UB}$, respectively.  Then  Galois theory implies that
 $ \E_{\VB}^*= [\E_{\VB}^*: 1 ]=  [\D: \E_{\VB}]$ is even. 
We conclude that the unique involution 
$\iota $ of $\Gamma$ is an element of $\E_{\VB}^*$.  Therefore,  $\iota$ fixes the field 
$\E_{\VB}$ pointwise.  So $\iota$ fixes, to within isomorphisms, each of the $\E N$-modules $\VB^j$.
Because $\iota$ acts non--trivially on $\D$ and fixes $\E_{\VB}$,  it  must act non--trivially on  $\E_{\UB}$. 
We conclude that $\iota$ viewed as an $\F$-automorphism of  $\E_{\UB}$ must coincide with the unique 
involution in  the Galois group $\Gal(\E_{\UB}  / \F)$  of $\E_{\UB}$ above $\F$. 
Furthermore, this unique involution must send $\UB ^i$ to its dual $\widehat {\UB ^i}$, for every $i=1, \dots, n$.
(Of course $\UB ^i $ is not self--dual, because it is a non--trivial absolutely irreducible 
module of the odd order group $G$ (see Lemma \ref{intr5}).) Hence,    
applying $\iota$ to both sides of \eqref{ee4} we get 
$$
\widehat{\UB^1}_N \cong e'(\widehat{\VB^1}  \oplus \dots \oplus \widehat{\VB^c}) \cong e'
(\VB^1 \oplus \dots \oplus \VB^c). 
$$
 Hence   the dual $\widehat{\VB^1}$ of $\VB ^1$
should be among the $G$-conjugates $\VB^1, \dots, \VB^c$ of $\VB^1$.  Because $\VB^1$ is not self--dual  
the $G$-conjugates of $\VB^1$ should appear in dual pairs. Hence  $c$ is even.
But $c$ is also  odd as a divisor of $|G|$.  This contradiction implies that $e$ is odd.
So the lemma holds.
\end{proof}

\begin{lemma}\mylabel{l7}
The group $G$ is not abelian.
\end{lemma}

\begin{proof}
Assume that $G$ is abelian. Then any cyclic  subgroup $N= <\sigma>$ of  $G$ is normal, for every $\sigma \in G$. 
Because $\B_N$ is hyperbolic, 
Lemmas \ref{l1} and  \ref{l3} along with  Corollary \ref{intr4} imply that 
\begin{equation*}\mylabel{e8}
\B_N  \cong 2 \cdot \Delta_{(N)},  
\end{equation*}
where $\Delta_{(N)} $ is a semi--simple $\F N$-submodule of $\B$. Using the splitting field $\E$ of $G$,  we write 
$\B  ^{\E}$ for the extended $\E G$-module $\B  ^{\E} = \B \otimes_{\F} \E$.  Then 
\begin{equation}\mylabel{e9}
\B ^{\E}_N \cong 2 \cdot \Delta_{(N)}^{\E}, 
\end{equation} 
where $\Delta_{(N)}^{\E}$ is the extended $\E N$-module $\Delta_{(N)} \otimes_ {\F} \E$. 
Let $\phi$ be a  Brauer character that the $\E G$-module $\B^{\E}$  affords. 
Because $(p, |G|) = 1$,  $\phi$  is defined for every element of $G$. So $\phi$ coincides with a complex character of $G$.
In view of \eqref{e9}, for every  cyclic subgroup $N= <\sigma>$ of $G$,   the restriction
 $\phi_N$ of $\phi$ to $N$ equals $2 \cdot \delta_{(N)}$,  where
$\delta_{(N)}$ is a complex character of $N$.
Hence, for every element $\sigma \in G$, the integer $2$ divides $\phi(\sigma)$  in the ring $\ZZ [\omega]$, 
where $\omega$ is a $|G|$-primitive root of unity. We conclude that $2$ also divides 
$\sum_{\sigma \in G} \phi(\sigma) \cdot \lambda(\sigma ^{-1})$, for any  irreducible  (linear) complex  character $\lambda$ of $G$.
That is, $2$ divides $|G| \cdot <\phi, \lambda>$, for any $\lambda \in \Irr(G)$. The fact that 
 $G$ has odd order, implies that  $2$ divides $<\phi, \lambda>$ in $\ZZ[\omega]$,  for any $\lambda \in \Irr(G)$.
Because $\phi = \sum_{\lambda \in \Irr(G)} <\phi, \lambda> \cdot \lambda $, we get 
\begin{equation}\mylabel{e12}
\phi = 2 \cdot \chi, 
\end{equation}
where $\chi$ is a complex character of $G$.

On the other hand, Lemma \ref{l3} implies that $\B = \UB_1 \oplus  \dots \oplus \UB_k$, where the 
$\UB_i$ are distinct simple $\F G$-modules, 
for all $i=1, \dots, k$. Hence  the extended $\E G$-module $\B ^{\E}$  will also
equal the direct sum of the distinct   $\E G$--modules \,$\UB_{1}^{\E} ,
\ldots ,\UB_{k}^{\E}$.  
 By Theorem 9.21 in \cite{3}, for each 
$ i=1, \ldots , k$, there exist
absolutely irreducible $\E G$--modules  $\Xii _i ^j$, for $j=1, \ldots ,
n_i $ such that 
$$
{\UB _i} ^{\E} \cong \bigoplus_{j=1}^{n_i} {\Xii_i}^j.
$$
Furthermore, the $\Xii_i ^j$, for $j=1, \dots, n_i$,  constitute a Galois
conjugacy class over $\F$, and thus they are all distinct.  In addition, the above 
absolutely irreducible  $\E G$-modules $\UB_i^j$,  for all $i=1, \dots, k$ and all $j=1, \dots, n_i$,  are distinct.
Indeed,   for all $i=1, \dots, k$,  the 
corresponding  simple $\F G$-modules $\UB_i$ are distinct.  
We conclude that  
$$
\B^{\E} \cong \bigoplus_{i=1}^k  \bigoplus_{j=1}^{n_i} {\Xii_i}^j,
$$
where $\UB_i^j$ are all distinct absolutely irreducible $\E G$-modules.
So the  character  $\phi$ that $\B ^{\E}$ affords  equals 
$$
\phi = \sum_{i=1}^k \sum_{j = 1}^{n_i} \chi_i^j, 
$$
 where,   for all $i=1, \dots, k$ and all $j=1, \dots, n_i$,  the character $\chi_i^j$ is a  Brauer character 
that $\UB_i^j$ affords.  So all these characters are distinct. This contradicts \eqref{e12}. 
Hence the group $G$ is not abelian, and the lemma is proved.
\end{proof}

\begin{lem}\mylabel{l4}
$G$ acts faithfully on $\mathcal{B}$.
\end{lem}

\begin{proof}
Suppose not. Let $K$ denote the kernel of the action of $G$ on
$\mathcal{B}$ and $\bar{G} = G / K$. Thus
 $ |\bar{G}| \lneq |G| $ (as $K\ne {1}$ ).

If $\bar G$ is not itself cyclic,  then any cyclic    subgroup $\bar N$
 of $\bar G$ is the image $\bar N = N/G$ of some proper subgroup $N$ of $G$.
Since $\B$ is $\F N $-hyperbolic, it is clearly $\F \bar{N}$-hyperbolic. 
 Hence
the triplet ${\mathcal{F, B}}, \bar{G}$ satisfies the hypothesis of
the Main Theorem. The minimality of $|G|$ implies that  $\mathcal{B}$ is
a hyperbolic  ${\mathcal{F}}\bar{G}$--module, and therefore a
hyperbolic  ${\mathcal{F}}G$--module, because any 
${\mathcal{F}}\bar{G}$--submodule of $\mathcal{B}$ is also an
${\mathcal{F}}G$--submodule of  $\mathcal{B}$. This contradicts the Inductive
 Hypothesis.

If $\bar G$ is cyclic,  then $\bar G = <\bar {\sigma}>$,  where $\bar{\sigma}$ is the image in $\bar{G}$ of 
 some $\sigma \in G$. Let 
$M = <\sigma>$.  Then $M$  is a proper subgroup of $G$,  
because $G$ is not cyclic.  In addition,  the image of  $M$ in $\bar G$ is $\bar G$. 
So $G =M K$ with $M \lneq G$. Then Remark  \ref{r1}  implies that  $\B$ is
$\F M$--hyperbolic and thus $\F G$--hyperbolic.
This last contradiction  implies the  lemma.
\end{proof}

\begin{lem}\mylabel{l5}
Suppose $M$ is a minimal normal subgroup of $G$. Then  $M$ 
is cyclic and central.
\end{lem}

\begin{proof}
According to Lemmas \ref{l3}  and \ref{l3.3}
 for each $ i=1, \ldots ,k $  there is a simple
${\mathcal{F}}M$--submodule ${\mathcal{V}}_i$ of  ${\mathcal{U}}_i$
and an odd  integer $e_i$, such that  ${{\mathcal{U}}_i}|_M
\cong e_i{\mathcal{V}}_i$. As $G$ acts faithfully on $\B$, there is some 
$i \in \{1, \dots, k\}$ such that $\VB_i \ne 1$ is non-trivial.
Let $K_M(\VB_i)$ be  the kernel of the action of $M$ on $\VB_i$.
The fact that $\VB_i$ is $G$-invariant implies that  
 $K_M(\VB_i)$ is a normal subgroup of $G$ contained in $M$. 
Hence $K_M(\VB_i)=1$. Therefore  $M$ admits a faithful simple representation. 
In addition, $M$ is a $q$-elementary abelian group, for some prime $q$ that 
divides $|G|$, because $G$ is solvable. 
We conclude that $M \cong \ZZ_q$ is a cyclic group of order $q$.

It remains to show that $M$ is central. If  
$\F$ is a splitting field of $M$
(that is, it contains a primitive $q$-root of $1$), 
then the fact that there exists a faithful, simple and thus one--dimensional, 
 $G$-invariant  $\F M$-module $\VB_i$  implies that $M$ is central in $G$.
If $F$ is not a splitting field of $M$,  we  work  with  the extension field $\E$
 of $\F$.
The extended module $\B ^{\E}= \B \otimes_{\F} \E $ equals  the direct sum of the extended $\E G$--modules \,$\UB_{1}^{\E} ,
\ldots ,\UB_{k}^{\E}$, because $\B$ is the direct sum of $\UB_1, \dots, \UB_k$. 
 As we have already seen,  for each 
$ i=1, \ldots , k$, there exist
absolutely irreducible $\E G$--modules  $\Xii _i ^j$, for $j=1, \ldots ,
n_i $,  that constitute a Galois conjugacy class over $\F$ and satisfy  
\begin{equation}\mylabel{e0}
{\UB _i} ^{\E} \cong \bigoplus_{j=1}^{n_i} {\Xii_i}^j.
\end{equation}
Since 
$\UB_i|_M \cong e_i \VB_i $ we have  $\UB_i^{\E}|_M  \cong e_i \VB_i^{\E}$.
 In addition, 
$$
\VB ^{\E}_i  \cong \bigoplus_{r=1}^{s_i} {\W}_i^r,  
$$
where the   ${\W}_i^r$, for  $r=1, \dots, s_i$, are absolutely irreducible  $\E M$-modules, and thus 
of dimension one,     that form a  Galois conjugacy class  over $\F$.
Therefore, 
\begin{equation} \mylabel{e1}
\UB_i^{\E}|_M \cong \bigoplus_{j=1}^{n_i}\Xii_i^j|_M \cong e_i \bigoplus_{r=1}^{s_i} {\W}_i^r, 
\end{equation}
for all $i=1, \dots, k$.

 As we have already seen, there exists
  $i \in \{ 1, \dots, k\} $ such that $\VB_{i}$ 
is a faithful $\F M$-module. Without loss,  we may assume that $i=1$. 
Then it is clear that the  ${\W}_{1}^r$ are  faithful $\E M$-modules, 
 for all $r \in \{1, \dots, s_1\}$. 
If ${\W}_{1}^r$  is $G$-invariant, for some $r \in \{1, \dots, s_1\}$ 
(and thus for all such $r$)  we are done.

Thus we may assume that the stabilizer $G_{\W}$ of ${\W}= {\W}_{1}^1$ in $G$ is 
strictly smaller  than $G$. 
Then $G_{\W}  = C_G(M)$, because ${\W}$ is 
$\E M$-faithful.
Let $C := G_{\W} = C_G(M)$.
 Note that $C $ is a normal subgroup of $G$, since 
$M \unlhd G$.  Furthermore, $C$  is also the stabilizer of ${\W}_1^r$, for all 
$r= 1, \dots , s_1$. According to   Lemma \ref{l3.3}, for all $i = 1, \dots, k$,   we have 
 $\UB_i|_C = m_i \cdot \XB_i  $,  where $\XB_i$  
is a simple $\F C$-module, and $m_i$ some positive integer.
For the extended $\E C$-modules $\XB_i^{\E}$  we  have 
$$
\UB_i^{\E}|_C \cong    m_i\XB_i^{\E} \cong m_i \bigoplus_{l=1}^{t_i} \XB_i^l,  
$$
where the $\XB_i^l$, for $ l=1, 2, \dots, t_i$,  are absolutely irreducible $\E C$-modules that constitute 
a Galois conjugacy class over $\F$. 
Hence 
\begin{equation}\mylabel{e2}
\UB_i^{\E}|_C \cong \bigoplus_{j=1}^{n_i}\Xii_i^j|_C \cong  m_i
 \bigoplus_{l=1}^{t_i} \XB_i^l, 
\end{equation}
for all $i=1, \dots, k$. We remark here  that,  because $\UB_i$ is quasi--primitive,  all the group 
conjugates of $\XB_i^1$ are among its Galois conjugates, for every $i=1, \dots, k$.

In the case $i=1$, equations \eqref{e1} and \eqref{e2}  imply 
\begin{equation}\mylabel{e3}
\UB_1^{\E}|_M  \cong   m_1\bigoplus_{l=1}^{t_1} \XB_1^l|_M   \cong 
e_1\bigoplus_{r=1}^{s_1} {\W}_1^r. 
\end{equation}
Without loss we may assume that $\Xii_1^1$ lies above $\XB_1^1$, and that  $\XB_1^1$ lies 
above ${\W}_1^1= {\W}$.  Clearly $\XB_1^1$ is non--trivial as it restricts to a multiple of the non--trivial $\F M$-module $\W_1$.
Hence Lemma \ref{l3.3} implies that $m_1 $ is an odd integer. 
Because $C$ is the stabilizer of $ {\W}$ in $G$, Clifford's
 theory implies that  $\XB_1^1$  is the unique simple  $\E C$-module  that lies above 
${\W_1^1}$ and induces irreducibly to  $\Xii_1^1$ in $G$.
Note that $\XB_1^1$ appears with odd multiplicity $m_1$ as a summand of $\UB_1^{\E}|_C$, 
because  the  $\E C $-modules $\XB_1^l$  are distinct for distinct values of $l$, 
as they form a  Galois  conjugacy class  over $\F $. 
Furthermore, if  $\XB_1^1$  lies under  some $\Xii_i^j$,  for $i  \ne 1$, then it 
 induces $\Xii_i^j$.  So $\Xii_i^j  \cong  \Xii_1^1$.  Hence the sum of the Galois 
conjugates of $\Xii_i ^j$  is isomorphic to the sum of the Galois conjugates of $\Xii_1^1$. 
Therefore
$$
\UB^{\E}_i \cong \bigoplus_{j=1}^{n_i}\Xii_i^j \cong
 \bigoplus_{j=1}^{n_1}\Xii_1^j  \cong \UB^{\E}_1. 
$$
The above contradicts the fact that $\UB_1$  and  $\UB_i$ are non-isomorphic
  simple $\F G$-modules (see Lemma \ref{l3}).  
We conclude that $\XB_1^1$ appears with odd multiplicity   $m_1$  in the decomposition of 
$$
\B^{\E}|_C  \cong \bigoplus_{i=1}^{k} \UB_i^{\E}|_C \cong  \bigoplus_{i=1}^k 
m_i \bigoplus_{l=1}^{t_i} \XB_i^l.
$$

On the other hand,  in view of  Corollary \ref{intr4} every 
simple $\F C$-submodule of 
$\B$ appears with even multiplicity in 
any decomposition of $\B_C$, as
$C$ is a normal subgroup of $G$ and  $\B_C$ is hyperbolic as an $\F C$-module, by Remark 
\ref{r1}.  Hence every absolutely irreducible $\E C$-submodule of $\B^{\E} $ should 
 also appear  with  even multiplicity in any decomposition of $\B^{\E}|_C$.
This  contradicts the conclusion of the preceding paragraph. So we must have  $G_{\W} = C = G$. 
Hence the lemma is proved.
\end{proof}

Clearly  Lemma \ref{l5} implies 
\begin{cor}\mylabel{co2}
Suppose that $M$ is a minimal normal subgroup of $G$ and $\E$ a  splitting  field of $G$ and all its subgroups. 
Then every $\E M$-module is $G$-invariant. 
\end{cor}

We can now show 
\begin{lemma}\mylabel{l5.5}
Suppose that $M $ is a minimal normal subgroup of $G$. Then the restriction 
$\B_M$ is homogeneous. Furthermore $\B_M \cong e \VB$, where $\VB$ is a simple  
faithful  $G$-invariant $\F M$-submodule of $\B_M$  and $e$ is a positive integer.
\end{lemma}

\begin{proof}
As in the previous lemma we write $\UB_i |_ M  = e_i \VB_i$, 
where $i= 1, \dots, k$, 
and $\VB_i$ is a simple $G$-invariant  $\F M$-submodule of $\UB_i$. 
If $\F M$ is not  homogeneous, then there  are at least two non-isomorphic 
simple $\F M$-submodules of $\B _M $, say $\VB$ and $\mathcal{W}$. 
We may  suppose that $\VB$ is non-trivial.
Assume that ${\mathcal{V}}_i \cong \mathcal{V}$ as ${\mathcal{F}}M$--modules, 
for all $ i=1, \ldots, l $ and  some $l$ such that $1\leq l\leq k$,  while 
${\mathcal{V}}_i \ncong \mathcal{V}$ for $i=l+1, \ldots ,k$. Let
${{\mathcal{U}}}$ be the  orthogonal direct sum 
$$
{\mathcal{U}}={\mathcal{U}}_1 \dot{\perp}
\ldots \dot{\perp}  \, {\mathcal{U}}_l.
$$ 
of the 
corresponding
${\mathcal{F}}G$--submodules of ${\mathcal{B}}$.
We also write  
$$
\RB = \UB_{l+1} \dot{\perp}
\ldots \dot{\perp}  \, {\mathcal{U}}_k, 
$$
for the  orthogonal direct sum of the remaining   simple $\F G$-submodules of $\B$.
Clearly
$$
\B = \UB \dot{\perp} \RB,   
$$
 while $\UB_M$ and $\RB_M$ have no  simple $\F M$-submodules in common.

 We will show
\begin{claim}\mylabel{c1}
 ${\mathcal{U}}$ is ${\mathcal{F}}N$-hyperbolic for every  cyclic subgroup $N$ of $G$.
\end{claim}

We first prove Claim  \ref{c1}  in the case that the  product
$N M$ is a proper subgroup of $G$. 
 In this case Remark~\ref{r1} implies that $\B _{NM}$ is
hyperbolic. Hence there exists  a self--perpendicular $\F NM$-submodule
$\SB >0$ of $\B$.
Then $\SB$ is  a maximal  isotropic $\F NM$-submodule of $\B_{NM}$.
Furthermore, 
 $\B_{NM} = \UB_{NM} \dot{\perp} \RB_{NM}$, where $\UB_{NM}$ and $\RB_{NM}$
have no simple $\F NM$-submodule in common  (otherwise $\UB_M$ and $\RB_M$ 
would have some common simple $\F M$-submodule).
Hence 
$$
\SB =  (\SB \cap \UB_{NM})  \dot{\perp} (\SB \cap \RB_{NM}).
$$
Because $\SB$ is isotropic, both  $\SB \cap \UB_{NM}$ and  $\SB \cap \RB_{NM}$ are also isotropic.
Hence their  $\F $-dimensions are at most $1/2$ the dimensions of $\UB_{NM}$ and 
$\RB_{NM}$,  respectively.  But $\SB$ is self--perpendicular and thus its $\F$-dimension 
is  exactly $(1/2) \dim (\B_{NM})$. We conclude that 
the $\F$-dimensions of $\SB \cap \UB_{NM}$ and  $\SB \cap \RB_{NM}$  are exactly 
$1/2$ the dimensions of $\UB_{NM}$ and 
$\RB_{NM}$,  respectively. 
 Therefore $\SB \cap \UB_{NM}$  is a maximal isotropic  $\F NM$-submodule of $\UB_{NM}$ of dimension $1/2$  the
 dimension of $\UB _{NM}$.  So $\SB \cap \UB _{NM}$ is
 self--perpendicular, by Lemma \ref{ll1}.  Thus $\UB _{NM}$ is hyperbolic as an
 $\F  NM$-module. Hence it is also hyperbolic as an $\F N$-module.
So Claim  \ref{c1} holds  when $NM  < G$.  
  
Assume now that $N$ is a cyclic subgroup of $G$ such that
 $NM = G$.  Because $M$ is minimal,  Lemma \ref{l5} implies  that
$M \cong  \ZZ_q$  is central. Hence  $G = M   N$ is an abelian group. 
This contradicts Lemma \ref{l7}. 
Therefore  $NM < G$,  for every  cyclic subgroup $N$ of $G$. Thus 
Claim \ref{c1} holds.  
 
Since   \, $\UB  < \B$,   the Inductive  Hypothesis,  along with 
Claim \ref{c1},  implies that $\UB$  is $\F G$-hyperbolic.
Hence $\UB$ contains a self--perpendicular  $\F G$-submodule $\TB$. 
Let $\TB^{\perp}$ be the submodule of $\B$ that is perpendicular to $\TB$.
Then $\RB $  as well as $\TB$ are subsets of  $\TB^{\perp}$. 
We conclude that $\TB$ is an isotropic  $\F G$-submodule of $\B$.
 Hence $\B$ is not anisotropic.    
This last contradiction implies that $\UB = \B$, and completes 
the proof of Lemma \ref{l5.5}.
\end{proof}

 \begin{lem} \mylabel{l6}
Every abelian normal subgroup of $G$ is cyclic.
\end{lem}

\begin{proof}
Let $A$ be an abelian normal subgroup of $G$.   
By Lemma \ref{l3.3} there is a simple $\F A$--submodule
${\mathcal{R}}_1$ of $ \UB _1$ and an integer $e_1$ such that 
$$
{\UB_{1}}|_{A} \cong e_1 {\mathcal{R}}_1.
$$
It follows from Lemma \ref{l5.5}  that $\mathcal{R}_1$ is non-trivial, since its restriction 
to any minimal  normal subgroup  of $G$ is non--trivial. 
Let $K_1$ denote the
corresponding  centralizer of ${\mathcal{R}}_1 $ in $A$ . Then $K_1 $
equals the centralizer $C_A(\UB _1)$ of $\UB _1$ in $A$,  and therefore is a
normal subgroup of $G$.
If $K_1$ is not trivial then it contains a minimal normal subgroup
 $M$ of $G$.  In view of Lemma~\ref{l5.5} the restriction ${ \UB _1}|_{M}$, cannot be
 trivial,  contradicting the definition of $K_1$. Hence $K_1$ is
 trivial. Thus $A$ is cyclic and the lemma is proved.
\end{proof}

Let $F = F(G)$ be the Fitting subgroup of $G$. Assume further that $\{ q_i\}_{i=1}^r$ are 
the distinct primes dividing $|F|$, and that  $T_i$ is the $q_i$-Sylow subgroup of $F$, 
for each $i=1, \dots, r$. Then  $F = T_1 \times T_2 \times \dots  \times T_r$.
Every characteristic abelian subgroup of $F$ is cyclic, according to 
Lemma \ref{l6}. Hence (see Theorem 4.9 in \cite{2})  either  $T_i$ is  cyclic 
or $T_i$ is the central product  $T_i = E_i \odot  Z(T_i)$
of the  extra special $q_i$--group
$E_i = \Omega (T_i)$  of exponent $q_i$
and the  cyclic group $Z(T_i)$.
We complete the proof of  Theorem A  exploring  the two possible types of $T_i$.

Assume  first that  $T_i$ is a  cyclic group, for all $i=1, \dots, r$.
In this case $F = T_1 \times \dots \times T_r$ is also a cyclic group.
Let $C/ F$ be a chief factor of $G$.  So $\bar{C} = C / F$ 
is an elementary abelian $q$-group, for some prime $q$, because $G$ is 
solvable. Then $\bar{C}$ acts  coprimely 
on $T_i$ for all   $i$  such  that  $q $ does not divide  $|T_i|$.
But $T_i$ is cyclic, and the minimal subgroup of $T_i$ is central in $G$. 
Hence $C_{T_i}(\bar{C}) \ne 1$. We conclude that
 $T_i = [T_i , \bar{C}] \times C_{T_i}(\bar{C}) = C_{T_i}(\bar{C})$. So  any  $q$-Sylow subgroup  $C_q$ 
of $C$ centralizes  the $q'$-Hall subgroup $R$ of $F$ that is also a $q'$-Hall subgroup of $C$.
We conclude that $C = C_q  \times R $. But $R$ is nilpotent as a subgroup of $F$.
 So $C$ is a nilpotent normal subgroup of $G$  bigger than  the Fitting subgroup $F$  of $G$.
Therefore  $G = F$ is a cyclic group,   contradicting the Inductive Hypothesis.
Hence  there exists  a Sylow subgroup  $T_i$ of $F = F(G)$  that is not cyclic. 

Let $T=T_i$ be a  non-cyclic  $q$-Sylow subgroup of $F$, where $q=q_i$ for some $i=1, \dots, r$.
Then  $T = E \odot Z(T)$,  where  $E = \Omega (T)$ is an extra special $q$-group
 of exponent $q$ and $Z(T)$ is the center of $T$. 
Of course $E$ is a normal subgroup of
$G$,  since   it is a characteristic subgroup of $F$.
Furthermore, $Z(E)$ is a central subgroup of $G$  because it is a minimal (it has order $q$)
normal  subgroup of $G$.
According to Lemma \ref{l5.5}, there exists a faithful $G$-invariant  $\F Z(E)$-module  $\VB$ 
so that   the restriction $\B _{Z(E)}$ of $\B $
to $Z(E)$  is a multiple of $\VB$.

Using the extension field $\E$  of $\F$,  
we  write $\VB ^{\E}$ for the extended $\E Z(E)$-module
$\VB \otimes_{\F} \E$.  Then 
\begin{equation*}
\VB^{\E} \cong \bigoplus_{j=1}^s  \VB^j
\end{equation*}
where $\VB^ j$ is an absolutely irreducible $\E Z(E)$-module, for all $j $ with 
$j= 1, \dots, s $.
Furthermore,  the $\VB^j$ constitute  a Galois conjugacy class over $\F$, and thus 
they are all distinct.  As  we have already seen (see  Corollary  \ref{co2} and Lemma \ref{l5.5}), the module 
 $\VB^j$ is a  non-trivial   $G$-invariant   $\E Z(E)$-module.
Because $E$ is extra special,    there exists a unique, up to isomorphism, absolutely irreducible $ \E E$-module 
 $\WB^j$     lying  above $\VB^j$, for every  $j=1, \dots, s$. 
  Note that  for all such $j$ 
the $\E E$-module  $\WB^j$  is $G$-invariant  because  $\VB^j$  is $G$-invariant.
According to Theorem  9.1 in \cite{4} (used for modules) there exists a canonical conjugacy class of subgroups $H \leq G$   
such that  $H E = G$  and $H \cap E = Z(E)$.  Furthermore,  for this conjugacy class 
there exists a one--to--one correspondence  between the isomorphism classes of   absolutely irreducible 
$\E G$-modules   lying above ${\WB}^j$ and those classes of  absolutely irreducible $\E H$-modules 
lying above $\VB^j$.  In addition, the fact that $G$ has odd order implies that   if $\Xi$ and $\Psi$ are representatives 
of the above two isomorphism classes,  then they correspond 
 iff $\Xi_H \cong  \Psi  \oplus 2 \cdot  \Delta$, 
where $\Delta$ is a  completely reducible $\E H$-submodule of $\Xi_H$.

Let $\UB = \UB_1 $  be one of the simple $\F G$-submodules of $\B$  appearing in \eqref{e-1}.
Then $\UB^{\E}  \cong   \oplus _{j=1}^{n_1} \Xii^j$, where the  $\Xii ^j$  are absolutely irreducible 
$\E G$-modules that form a Galois conjugacy  class.  
 As earlier, we write  $\E_{\UB}$ for the extension field  of $\F$  generated by 
all the values of the absolutely irreducible character that $ \Xii ^1$ affords.
Let  $\Gamma = \Gal( \E_{\UB}/ \F)$   be the  Galois group of that extension. 
Then (see Theorem 9.21 in \cite{3}), 
\begin{equation}\mylabel{e4}
 \UB^{\E}  \cong   \bigoplus_{j=1}^{n_1} \Xii ^j  \cong  \bigoplus _{\tau \in \Gamma} (\Xii^1) ^{\tau}, 
 \end{equation}
Clearly $\Xii ^1 $  lies above  ${\WB}^j$, for some  $j=1, \dots, s$, since $\UB  = \UB_1$ lies above $\VB$. 
Let $\Psi$ be a representative of the isomorphism class of absolutely irreducible $\E H$-modules 
that  corresponds to $\Xii ^1$ and lies above $\VB^j$.
Then 
\begin{equation}\mylabel{e5}
\Xii^1_H  \cong \Psi   \oplus  2 \cdot \Delta,
\end{equation}
for some completely reducible $\E H$-module $\Delta$.
Let $\E_{\Psi}$   be  the subfield of $\E$ generated  by $\F$ and  all the values 
of the absolutely irreducible character that $\Psi$ affords. Then $\E_{\Psi}$ is a Galois extension of $\F$.
Furthermore, 
 \begin{equation}\mylabel{e6}
\E_{\Psi} = \E_{\UB}.
\end{equation}
Indeed,   for any element $\sigma $ in the Galois group $\Gal(\E / \F)$  of $\E$ above $\F$ 
 we get 
$$
(\Xii^1)^{\sigma}_H \cong  \Psi^{\sigma}   \oplus 2 \cdot \Delta^{\sigma}.
$$
Hence   $(\Xii^1)^{\sigma}$ corresponds  to $\Psi^{\sigma}$, as $\Psi^{\sigma}$  is the  only absolutely irreducible 
  $\E H$-module that appears with odd multiplicity in $(\Xii^1)^{\sigma}_H$.
Therefore,  $(\Xii^1)^{\sigma} \ncong \Xii^1 $ iff  $\Psi^{\sigma} \ncong \Psi$. 
This is enough to guarantee  that  \eqref{e6} holds.
We conclude that the  sum $\oplus_{\tau \in \Gamma} \Psi^{\tau}$  
 is the extension to $\E $ of  an irreducible $\F H$-module,  i.e.,  there exists an irreducible 
$\F H$-module  $\Pi$ such that 
$$
\Pi^{\E}  \cong \bigoplus_{\tau \in \Gamma} \Psi^{\tau},   
$$
where $\Pi^{\E}$ is the extended  $\E H$-module $\Pi \otimes_{\F} \E$.
Furthermore, \eqref{e4} and \eqref{e5} imply that  
 $\Pi$ appears with odd multiplicity as a summand 
 of  $\UB_H = \UB_1|_{H}$.  

Next we observe that  if 
 $\Pi $ appears as a summand  of   $\UB_i |_{H}$, for some $i=2, \dots, k$, 
then it appears with even multiplicity.  The reason is that
 $\UB_1  \ncong \UB_i$ for all such  $i$.
As in  \eqref{e0} we choose a Galois conjugacy class 
$\{\Xii_i^j\}_{j+1}^{n_i}$ of absolutely irreducible $\E G$-modules  
such that $\UB_i^{\E} \cong \oplus_{j=1}^{n_i} \Xii_i^j$.  Then 
 $\UB_i  \ncong \UB  = \UB_1$ implies that 
 $\Xii_i ^j \ncong \Xii^1$, for all $i=2, \dots, k$ and all $j=1, \dots, n_i$.
So  the $\E H$-module  $\Psi$ can't correspond to $\Xii_i^j$, for any such $i, j$. 
Therefore  if $\Psi$  appears  as a  summand of the restriction $\Xii_i^j|_{H}$ of 
$\Xii_i^j $ to $H$, then it appears only with even multiplicity. 
Hence the same holds for $\Pi$, i.e.,   $\Pi$ 
appears only with even multiplicity   as a summand of $\UB_i|_H$, 
whenever $i =2, \dots, k$. 
We conclude that  $\Pi$ appears with odd multiplicity as a summand of $\B_H  = \UB_1|_H \oplus \dots \oplus  \UB_k|_H$.

We complete the proof of  Theorem A with one more contradiction, that follows the fact that   
 $\Pi$ is a self--dual $\F H$-module.  That we get a contradiction if $\Pi$ is self--dual is easy to see, 
because according to   Proposition \ref{intr3}, $\Pi$ should appear with even multiplicity as a summand of
the hyperbolic $\F H$-module  $\B _H$.
Thus it suffices to show that $\Pi$ is self--dual.

The fact that $\UB= \UB_1$ is self--dual implies that $\UB^{\E}$ is also self--dual.
Hence the dual $\widehat{\Xii^1}$  of $\Xii^1 $  is a Galois conjugate $(\Xii^1) ^{\tau}$ to  $\Xii^1$,  for some $\tau  \in \Gamma$.
Furthermore,  \eqref{e5} implies that 
$$
\widehat{\Xii^1}_H \cong \widehat{\Psi}   \oplus  2 \cdot \widehat{\Delta}.
$$
Thus the dual $\widehat{\Xii^1}$  corresponds to the dual $\widehat{\Psi}$ of $\Psi$.
Therefore the dual $\widehat{\Psi}$ of $\Psi$  is a Galois conjugate of $\Psi$. 
Hence $ \Pi^{\E} \cong \oplus_{\tau \in \Gamma} \Psi^{\tau}$ is a self--dual $\E H$-module.
So $\Pi$ is also  self--dual.

This completes the proof of Theorem A.

\end{document}